\documentclass[12pt,a4paper]{article}

\usepackage{amssymb,amsmath, amsfonts, stmaryrd, amsthm, latexsym}
\usepackage{graphicx}
\usepackage[all]{xy}
\usepackage{pstricks}
\usepackage[latin1]{inputenc}
\usepackage{setspace}

\DeclareMathOperator{\Sub}{Sub}

\theoremstyle{plain} 
\newtheorem{satz}{Satz}[section]
\newtheorem{prop}[satz]{Proposition}
\newtheorem{lem}[satz]{Lemma}

\newtheorem{cor}[satz]{Corollary}
\newtheorem{thm}[satz]{Theorem}

\theoremstyle{definition}
\newtheorem{defi}[satz]{Definition}

\newtheorem{obs}{Observation}[section]

\newtheorem{rem}{Remark}[section]

\newcommand{\Ra}{\Rightarrow}

\newcommand{\0}{\{0\}}

\newcommand{\btl}{\blacktriangleleft}
\newcommand{\btr}{\blacktriangleright}
\newcommand{\gdw}{\Leftrightarrow}
\newcommand{\ra}{\rightarrow}

\newcommand{\hra}{\hookrightarrow}
\newcommand{\thra}{\twoheadrightarrow}

\newcommand{\eps}{\epsilon}
\newcommand{\vareps}{\varepsilon}

\newcommand{\sub}{\subseteq}

\newcommand{\ran}{\rangle}
\newcommand{\lan}{\langle}

\newcommand{\LR}{\overline{L}(R_R)}
\newcommand{\RL}{\overline{L}(_R R)}

\newcommand{\LS}{\overline{L}(S_S)}

\newenvironment{bew}{\noindent $\btr$ {\bf Proof. }}{$\btl$}
\newenvironment{bew2}{\noindent $\triangleright$ {\sf Proof.
  }}{$\triangleleft$} 
\newenvironment{abs}{\begin{center} {\sc Abstract.} \end{center} }

\newcommand{\lab}[1]{\label{#1}}

\begin{document}

\begin{center}
{\Large On linear representation of $\ast$-regular rings having
  representable
ortholattice of projections}\\[3mm]
{\large Christian Herrmann and Niklas Niemann}

\end{center}

\begin{abs}

The origings of regular and $*$-regular rings lie in the
  works of J. von Neumann and F.J. Murray on operator algebras, von
  Neumann-algebras and projection lattices. They
  constitute a strong connection between operator theory, ring theory
  and lattice theory.\\

This paper aims at the following results:
\begin{enumerate}
\item The class of all $*$-regular rings forms a variety.
\item A  subdirectly irreducible
  $*$-regular ring $R$ 
 is faithfully representable 
 (i.e. isomorphic to a subring of an endomorphisms ring of
  vector spaces, where the involution  is given by
  adjunction with respect to a scalar product on the vector space)
if so is its ortholattice of projections. 
\end{enumerate}
This is a short version of the second author's PhD  thesis. 
\end{abs}

\tableofcontents
\section{Introduction}\lab{Preliminaries}

In the present section, we will fix notational conventions
and introduce the different concepts of representability.\\

The term {\sl ring} is used for rings with or without unit. Rings are
denoted by $R,S,T,C$. We consider rings with an involution $^*:R \ra
R$ (an anti-automorphism of order two). We consider the unary map as
part of the signature of the ring $R$.  
A (von Neumann-)regular ring is a ring such that every element $x$ has
at least one quasi-inverse $y$, i.e., for $x \in R$ there exists an $y \in R$
such that $xyx=x$. A $*$-regular ring is a regular involutive ring satisfying the implication $xx^*=0
\Ra x = 0$.

The term {\sl idempotent} is used for a ring
element $x$ satisfying $x^2=x$, the term {\sl projection} is used for a ring
element $x$ satisfying $x^2=x^*=x$. We use the letters $p,q$ for
projections and $e,f,g$ for idempotents and projections.\\

The term {\sl lattice} is used for a partially
ordered set with binary operations join and meet. These operations are
denoted by $+$ and $\cdot$. All lattices considered have a smallest
element $0$. By a {\sl bounded lattice}, we mean a
lattice with top and bottom. We use the terms {\sl interval} and {\sl
  section} of a lattice in the usual way. Intervals and sections are
bounded lattices in their own right, with the inherited operations. We
use the notation $a \oplus b$ or $\bigoplus a_i$ for the join of
independent elements $a$ and $b$ or for the join of the independent
family $\{a_i: i \in I\}$. We define the {\sl height} $h(L)$ of a
lattice as usual to be the supremum of all cardinalities $|C|-1$, $C$
a chain in $L$.  

In this paper, we deal mainly with {\sl modular} lattices. Of
particular interest are ({\sl relatively or sectionally}) {\sl
  complemented modular lattices} and ({\sl sectional}) {\sl modular
  ortholattices}. We use the abbreviations CML and MOL,
respectively. We denote the {\sl orthocomplementation} on a MOL $L$ by
$^\perp: L \ra L$.  

\newpage
\subsection{Rings and Lattices}

In this section, we recall well-known results about regular rings and
the connections between regular rings and complemented modular lattices.

\begin{thm}A ring with unit is regular if and only if the set of all
  its principal right ideals is a complemented modular lattice. 

If $R$ does not contain a unit, the equivalence holds for {\em
  complemented} replaced by {\em relatively complemented}. 
\end{thm}

For a ring $R$, we denote the set of
  all its principal right (left) ideals by $\LR$ (by $\RL$).

\begin{lem}\lab{lLR} A ring with unit (without unit) is $*$-regular
  if and only if $\LR$ is a (sectional) MOL.
\end{lem}

\begin{bew2} Folklore. If $R$ is $*$-regular, every principal right
  ideal is generated by a projection. The orthogonality on $\LR$ is given by
$$ aR \perp bR \gdw b^*a = 0$$
If $R$ contains a unit, then the orthocomplement of $eR$, $e$ a
projection in $R$, is given by $(1-e)R$. 
\end{bew2}

\begin{prop} If $R$ is regular, then the lattices $\LR$ and $\RL$ of
  principal right ideals and principal left ideals respectively, are anti-isomorphic. 

If $R$ is $*$-regular, $\LR$ and $\RL$ are isomorphic.
\end{prop}

\begin{bew2}See \cite{flor}, \cite{skor} and \cite{mae58}.
\end{bew2}

\begin{lem}\lab{lRsimpleimpliesLRsimple}If $R$ is a simple
  $*$-regular ring (without unit), then $\LR$  is a simple (sectional) MOL. 
\end{lem}

\begin{bew2} See \cite[Theorem 2.5]{proat}.
\end{bew2}

\begin{lem}\lab{lsimpleMOLsimpleintervals} In a simple MOL,
  each non-trivial interval $[0,a]$ is simple.   
\end{lem}

\begin{bew2}\cite[Lemma 2.2]{JonsCMLRep}.
\end{bew2}

\begin{lem}\lab{leRe*regularsubring} If $R$ is $*$-regular and $e$ a
  projection in $R$, then the set $eRe$ is a $*$-regular subring (with
  unit $e$) of $R$.  
\end{lem}

\begin{bew2} Since $e$ is a projection, $eRe$ is a subring and the
  involution  on $R$ restricts to an involution on $eRe$. For
  regularity, consider $x \in eRe$. Take an quasi-inverse $y$ of $x$ in
  $R$ and reflect that $eye$ is also an quasi-inverse of $x$. 
\end{bew2}

Let $R$ be a $*$-regular ring and $e$ a projection in $R$. We write
$R_e$ for the $*$-regular subring $eRe$. Furthermore, we 
define the  {\sl height} $h(R)$ by the height $h(\LR)$ of its
principal ideal lattice $\LR$.

\begin{lem}\lab{lJonIntervaleRe} Let $R$ be a $*$-regular ring and $e$ a projection in
  $R$. Then the lattice $\overline{L}((R_e)_{R_e})$ of all principal
  right ideals in $R_e = eRe$ is isomorphic to the section $[0,eR] \sub
  \LR$.
\end{lem}

\begin{bew2}\cite[Lemma 8.2]{JonsCMLRep}.
\end{bew2}

\begin{lem}\lab{leResimple} If $R$ is a simple $*$-regular ring and
  $e$ a projection in $R$, then $R_e$ is a simple $*$-regular ring
  with unit $e$. 
\end{lem}

\begin{bew2} It is left to show simplicity. For a non-vanishing ideal $A$ in $eRe$,
  consider the ideal generated by $A$ in $R$.
\end{bew2}

\begin{lem}\lab{lprojx}Let $R$ be a $*$-regular ring. Then for each
  $x \in R$, there exists a projection $e_x \in R$ such that $e_x x e_x = x$.
\end{lem}

\begin{bew2}Let $p_x$ be the projection that generates the left ideal
  generated by $x$ and $q_x$ be the projection that generates the right
  ideal generated by $x$. Take $e_x:= p_x \vee q_x$ to be the supremum
  in the lattice of all projections of $R$.
\end{bew2}

\subsection{Frames}

In this section, we recall the notion of perspectivity of elements of
a lattice and the concept of a frame. The reader familiar with frames of modular lattices might
  give the following descriptions only a short glance and then jump to
  Definition \ref{dstableorthogonalnkframe} of a {\sl stable 
  orthogonal frame} and  Corollary
  \ref{csimpleMOLstableorthogonalnkframe2} that a simple MOL of
  height at least $n$ contains a stable orthogonal $(n,k)$-frame.\\

Two elements $a,b$ of a lattice $L$ are called {\sl
  perspective (to each other)} if they have a common complement $c$ in
$L$. If $a$ and $b$ are perspective, we write $a \sim b$. We say that
$a$ is {\sl subperspective to $b$} or {\sl perspective to a part of
  $b$} if there exists an element $d \leq b$ such that $a \sim d$. We
write $a \lesssim b$. An element $c$ establishing a (sub)perspectivity
between elements $a$ and $b$ is called an {\sl axis of
(sub)perspectivity} between $a$ and $b$. If $a \lesssim b$, the part
  $d \leq b$ such that $a \sim d$ is called the {\sl image of $a$
  under the perspectivity between $a$ and $b$}. 
 
There exist different notions of a {\sl frame}. In
\cite{JvN}, von Neumann defined a {\sl homogeneous basis} for a CML
$L$ (p. 93) and a {\sl (normalised) system of axes of perspectivity}
for a given homogeneous basis (p. 118). The combined system was called a {\sl
  (normalised) frame} for $L$.  Equivalently, G. Bergmann and A. Huhn
introduced the notion of a {\sl n-frame} (originally, a {\sl
  $(n-1)$-diamond}) in a modular lattice (see the survey articles
\cite{DayGeo}, \cite{DayOrder} or the article of C. Herrmann in memory
of A. Day \cite{HerrInMemoriamDay}).  

The notion of a frame was subject to further development and
  generalisation. See \cite{JonsCMLRep} for the introduction of a {\sl
  partial frame}, a {\sl large partial frame} and a {\sl global frame}.
In \cite{JonsCMLRep}, J\'{o}nsson defined a {\sl large partial
  $n$-frame} in a bounded modular lattice $B$ to be a subset of $B$
  consisting of independent elements $a_0, \dots, a_{n-1}, d$ and the
  entries of a symmetric matrix $c=(c_{ij})_{i,j <n}$ such that  
the supremum of $a_0, \dots, a_{n-1}$ and $d$ equals the unit element
$1_B$, $d$ consists of a sum of finitely many elements each
of which is subperspective to $a_0$, and $c_{ij}$ is an axis of
perspectivity between $a_j$ and $a_i$.  

We adapt the definition of J\'{o}nsson in the following way: Decomposing
$d$ into $k$ summands, each of which is subperspective to $a_0$, we
incorporate these summands and their axes of subperspectivity to $a_0$
into the
frame. Furthermore, we demand that the spanning elements of the frame
are independent. 

\begin{defi}\lab{dlargepartialnkframe}{\bf Large partial $(n,k)$-frame}\\
A {\sl large partial frame of format $(n,k)$} in a bounded modular
lattice $L$ is a subset  
$$\Phi:=\{a_i, a_{0i}: 0 \leq i < n+k\} \sub L$$
such that the following conditions are satisfied.
\begin{enumerate}
\item $\bigoplus\limits_{i < n+k} a_i = 1_L$
\item $a_0 + a_i = a_0 \oplus  a_{0i} = a_i \oplus  a_{0i}$ for $i =
  1, \dots, n-1$ 
\item $a_0(a_i + a_{0i}) + a_i = a_i \oplus a_{0i} = a_0(a_i + a_{0i})
  \oplus a_{0i}$ for $i=n, \dots, n+k-1$
\end{enumerate}
That is,  $\Phi$ contains $n+k$ independent elements $a_i$ spanning
the lattice $L$ (condition (1)). Conditions (2) and (3) state
that $a_1, \dots , a_{n-1}$ are perspective to $a_0$ and $a_n, \dots ,
a_{n+k-1}$ are subperspective to $a_0$, where the axes of
(sub)perspectivity are just the $a_{0i}$. In particular, we have 
$a_0 \cdot a_{0i} = a_i \cdot a_{0i} = 0$ for all $i$. 

The axes of perspectivity between $a_i,a_j$ for indices $i,j <n$ can
be constructed via the axes $a_{0i}, a_{0j}$: We have $a_{ji} = 
[a_{0j} + a_{0i}] \cdot [a_j + a_i]$ and consequently, we have 
$a_{ki} = [a_{kj} + a_{ji}] \cdot [a_k + a_i]$ for $i,j,k <
n$. Likewise, we can construct the axis of subperspectivity $a_{ji}$
between $a_i$ and $a_j$ for indices $i,j$ such that $j < n$ and $n
\leq i < n+k$.

For short, we call $\Phi$ a {\sl large partial $(n,k)$-frame} or an
{\sl $(n,k)$-frame}, dropping the attribute {\sl large partial} for the ease
of notation and to avoid confusion with the notion of a
large partial $n$-frame in the sense of J\'{o}nsson.
\end{defi}

In the following, we state some helpful results and develop the
appropriate notion of a frame for a modular ortholattice.

\begin{lem}\lab{lCMLsubperspectivity} Let $L$ be a CML 
  and assume that $a_0,a,b \in L$ are elements such that $a_0 \leq a$,
  $a_0 \cdot b = 0$, and $b$ is subperspective to $a_0$. Then the relative complement $d$ of $ab$ in $[0,b]$ is
  subperspective to $a_0$ and $a \oplus d = a + b$.
\end{lem}

\begin{lem}\lab{lCMLframes}Let $L$ be a CML. If $L$ contains a large
  partial $n$-frame in the sense of J\'{o}nsson, then $L$ contains an
  $(n,k)$-frame. 
\end{lem}

\begin{bew2} Construct summands $a_i, n   \leq i < n+k$ of $d$ with
  the desired properties (subperspective to $a_0$ and such that $a_0, \dots ,
  a_{n+k-1}$ are independent) inductively. 
\end{bew2}

Next, we introduce the concept of a stable frame. The main difference
is that we incorporate all the axes of (sub)perspectivity (see Definition
\ref{dlargepartialnkframe}) and a
set of relative complements.

\begin{defi}\lab{dstablenkframe}{\bf Stable $(n,k)$-frame}\\
Let $L$ be a CML. A subset
$$\Phi = \{a_i, a_{ij}:  0 \leq i < n, 0 \leq j < n+k\} \cup \{z_{ij}:
j < n, n \leq i < n+k\} \sub L$$
will be called a {\sl stable $(n,k)$-frame $\Phi$ in $L$} if 
\begin{enumerate}
\item $\{a_i, a_{0i}: 0 \leq i < n+k\} \sub L$ is a $(n,k)$-frame in
 $L$
\item for $i,j \in I$, $i < n$,  $a_{ij}$ is the
  axis of (sub)perspectivity between $a_j$ and $a_i$
\item for each pair $(i,j)$ of indices with $j < n$ and $n \leq i <n+k$, the
element $z_{ij}$ is a complement of $b_{ji}$ in $[0,a_j]$, where $b_{ji}$ is
the image of $a_i$ under the subperspectivity $a_{ji}$ between $a_i$ and $a_j$.
\end{enumerate}
\end{defi}

\begin{lem}\lab{lCMLstableframe}Let $L$ be a CML.

If $L$ contains an $(n,k)$-frame, $L$ contains a stable $(n,k)$-frame.
\end{lem}

\begin{bew2}Choose the necessary relative complements.
\end{bew2}

\begin{defi}\lab{dorthogonalnkframe}{\bf Orthogonal $(n,k)$-frame}\\
Let $L$ be a MOL.\\
An $(n,k)$-frame $\Phi$ in $L$ is called an {\sl orthogonal
$(n,k)$-frame} if the following additional condition is satisfied: 
$$\forall i \in \{0, \dots, n+k-1\}. \qquad a_i^\perp = \sum\limits_{j
  \neq i} a_j$$
\end{defi}

\begin{defi}\lab{dstableorthogonalnkframe}{\bf Stable orthogonal
  $(n,k)$-frame}\\ 
Let $L$ be a MOL. A {\sl stable orthogonal $(n,k)$-frame} is a stable
frame such that
\begin{enumerate}
\item $\Phi$ as a frame satisfies the condition of Definition
  \ref{dorthogonalnkframe}, and 
\item the relative complements $z_{ij}$ are relative orthocomplements.
\end{enumerate}
\end{defi}

\subsubsection{Orthogonalisation of a $(n,k)$-Frame via
  J\'{o}nsson}\lab{ssorthogonalizationofankframe} 

Now we want to show that the notion of a (stable) orthogonal
$(n,k)$-frame is the appropriate one for a MOL: We will proove that a
given frame can be orthogonalised. We will base the proof on arguments
and results presented in \cite{JonsCMLRep}. In fact, one could choose
an alternative approach via ideas of Fred Wehrung, presented in
\cite{DimL}, using the notion of a {\sl normal equivalence} in a
modular lattice and the concept of a {\sl normal} modular lattice. 

\begin{lem}\lab{lMOLprojectiveelements} Let $L$ be a MOL and $a,b$
  projective elements in $L$. Then there exist four elements $b_0,
  b_1, b_2, b_3$ in $L$ such that $b$ is the direct orthogonal sum of
  $b_0, b_1, b_2, b_3$ and each $b_i$ is perspective to a part of $a$.
\end{lem}

\begin{bew2}The lines of argument follow J\'{o}nsson's proof of Lemma 1.4
  in \cite{JonsCMLRep}. The only difference is that we choose the 
  relative complements in J\'{o}nsson's proof to be relative
  orthocomplements in the considered intervals.
\end{bew2} 

\begin{lem}\lab{lMOLperspectiveelements} Let $L$ be a MOL and $a_0, a,
  b$ elements in $L$ such that $a_0 \leq a$, $a \cdot b = 0$, and $b
  \lesssim a_0$.  Then $b$ decomposes into a direct 
  orthogonal sum of five elements $b_0, b_1, b_2, b_3, b_4$ such that
  each $b_i$ is subperspective to $a_0$. 
\end{lem}

\begin{bew2} We choose $c:=b \cdot a^\perp$
  and $d$ as the relative orthocomplement of $c$ in $[0, (a+b)\cdot
  a^\perp]$. As part of $b$, $c$ is subperspective to
  $a_0$. Furthermore, one can show that $d$ is projective to a part $x
  \leq a_0$. Consequently, by Lemma \ref{lMOLprojectiveelements}, we
  can decompose $d$ into 
  the direct orthogonal sum of 4 elements $b_1, \dots, b_4$, each of
  which is subperspective to $a_0$.  Together with $b_0:=c$, we have
  the desired result. 
\end{bew2}

\begin{lem}\lab{lMOLnontrivialindependentelements} Let $L$ be a
  simple MOL and $a,b \in L$ non-trivial independent elements. Then
  there exist non-trivial elements $a_0, b_0$ with $a_0 \leq a, b_0
  \leq b$ such that $a_0$ and $b_0$ are perspective to each other. In
  particular, this holds if $b \leq a^\perp$. 
\end{lem}

\begin{bew2}Since $L$ is simple, the neutral ideal 
  generated by $a$ is the whole lattice. Then
  $b$  is the sum of finitely many elements,
  each of which is perspective to a part of $a$. Choose one such
  non-trivial summand as $b_0$ and the corresponding perspective part
  of $a$ as $a_0$.
\end{bew2}

\begin{lem}\lab{lsimpleMOLlpnframe} Let $L$ be a simple MOL with
  $h(L)\geq n$. Then there
  exists a large partial $n$-frame $\Phi$ (in the sense of J\'{o}nsson)
  such that the first $n$ elements $a_0, \dots, a_{n-1}$ of $\Phi$ are
  orthogonal, that is, we have 
$$a_k \leq (\bigoplus\limits_{i <k} a_i )^\perp$$
\end{lem}

\begin{bew2} By induction.
\end{bew2}

\begin{lem}\lab{lMOLorthogonalnkframe} Let $L$ be a MOL containing a
  $(n,k)$-frame $\Phi$ such that $a_0, \dots a_{n-1}$ are orthogonal, 
  that is, if for all $k < n$, we have $a_k \leq (\oplus_{i < k}
  a_i)^\perp$.  Then $L$ contains an orthogonal $(n,k')$-frame for
  some $k'$.
\end{lem}

\begin{bew2} By induction over the elements $a_n, \dots , a_{n+k-1}$
  and Lemma \ref{lMOLprojectiveelements}.  
\end{bew2}

\begin{cor}\lab{csimpleMOLstableorthogonalnkframe2} Let $L$ be a
  simple MOL of height at least $n$.  Then $L$ contains a stable
  orthogonal $(n,k)$-frame.  
\end{cor}

\subsubsection{Rings and Frames}

Combining Corollary \ref{csimpleMOLstableorthogonalnkframe2} with
Lemma \ref{lRsimpleimpliesLRsimple} and Lemma \ref{leResimple}, we get
the following results. 

\begin{cor}\lab{cRsimpleimpliesLRcontainssonkframe} Let $R$ be a
  simple $*$-regular ring with unit and $h(R) \geq n$. Then the
  MOL $\LR$ contains a stable orthogonal $(n,k)$-frame.  
\end{cor}

\begin{cor}\lab{cRsimpleimplieseRecontainsframe} If $R$ is a simple
  $*$-regular ring and $e$ a projection in $R$, then the MOL
  $\overline{L}({R_e}_{R_e})$ of principal right ideals of $R_e$ 
  contains a stable orthogonal frame.  
\end{cor}

\begin{rem} Clearly, the format of the frame depends on the height
  of $R_e$. 
\end{rem}

\subsubsection{Projectivity of Frames}\lab{ssprojectivityofframes}

It is well-known that global frames are projective. In this section,
we state similar results for the above introduced frames. We begin
with large partial $(n,k)$-frames. 

\begin{lem}\lab{lprojectivitynkframes} Let $K,L$ be CMLs, $f: K
  \twoheadrightarrow L$ a surjective $0$-$1$-lattice homomorphism and
  $\Phi \sub L$ a large partial $(n,k)$-frame in $L$. Then there
  exists a section $M \leq K$ and a set $\Psi \sub M$ such that 
\begin{enumerate}
\item $f_{|_M}: M \ra L$ is a surjective lattice homomorphism,
\item $\Psi$ is a large partial frame in $M$ of the same
  format as $\Phi$, 
\item $f[\Psi] = \Phi$.
\end{enumerate}
\end{lem}

\begin{bew2} Inductive process and appropriate choices of preimages.
\end{bew2}\\

Similarly, we have the following. 

\begin{lem}\lab{lprojectivitystablenkframes} Let $K,L$ be CMLs, $f: K
  \twoheadrightarrow L$ a surjective $0$-$1$-lattice homomorphism and
  $\Phi \sub L$ a stable $(n,k)$-frame in $L$. Then there exists a
  section $M \leq K$ and a set $\Psi \sub M$ such that
\begin{enumerate}
\item $f_{|_M}: M \ra L$ is a surjective lattice homomorphism,
\item $\Psi$ is a stable frame in $M$ of the same
  format as $\Phi$, 
\item $f[\Psi] = \Phi$.
\end{enumerate}
\end{lem}

\begin{bew2} Incorporate the choice of the necessary relative
  complements in the inductive procedure. To accomplish this, it is
  enough to show that if $a,b$ are in $K$ such that $b \leq a$ and
  $f(b) \oplus c   = f(a)$ for some $c \in L$, then there exists $d
  \in K$ such that $b \oplus d = a$ and $f(d) = c$.
\end{bew2}

\begin{lem}\lab{lprojectivitystableorthogonalnkframes} Let $K,L$ be
  MOLs, $f: K \twoheadrightarrow L$ a surjective $0$-$1$-lattice
  homomorphism and  $\Phi \sub L$ a stable orthogonal $(n,k)$-frame in
  $L$. Then there exists a section $M \leq K$ and a set $\Psi \sub M$ such
  that
\begin{enumerate}
\item $f_{|_M}: M \ra L$ is a surjective lattice homomorphism,
\item $\Psi$ is a stable orthogonal frame in $M$ of the same
  format as $\Phi$, 
\item $f[\Psi] = \Phi$.
\end{enumerate}
\end{lem}

\begin{bew2}Incorporate the orthogonality into the inductive process. 
\end{bew2}

\subsection{Concepts of Representability}

\begin{defi}\lab{dlinearrep}{\bf Linear representation}\\
As in  \cite{flor} and \cite{dip}, a {\sl linear positive representation} of a $*$-regular ring
$R$ is a tuple
$$\sigma = (D, V_D, \phi, \rho)$$ 
where $D$ is an (involutive) skew field, $V_D$ a right vector space
over $D$, $\phi$ a scalar product on $V_D$, and 
$$ \rho: R \ra End(V_D)$$ 
a ring homomorphism such that
$$ \forall r \in R. \qquad \rho(r^*) = \rho(r)^{*_{\phi}}$$

If the morphism $\rho: R \ra End(V_D)$ is injective, we call $\sigma$ a
{\sl faithful} representation.
\end{defi}

\begin{defi}\lab{dgrep}{\bf Generalised representation}\\
Let $R$ be an involutive ring, $I$ an arbitrary non-empty index set
and $\sigma$ a tuple 
$$\sigma = (I,\{D_i\}_{i \in I},\{V_i\}_{i \in I},\{\phi_i\}_{i \in
  I},\rho)$$
consisting of an indexed family of (involutive) skew
  fields, an indexed family of vector spaces and an indexed family of
  scalar products, such that for each $i \in I$, 
$V_i$ is a right vector space over $D_i$ with scalar product $\phi_i$, 
and a map 
$$\rho: R \ra \prod\limits_{i \in I} End(V_{iD_i}).$$

If $\rho$ is a $*$-ring morphism, i.e., for all $r \in R$ and all
$i \in I$ the condition
$$\pi_i(\rho(r^*)) =   \big(\pi_i(\rho(r))\big)^{*_{\phi_i}}$$
holds, we call $\sigma$ a {\sl positive generalised representation} of
$R$. For short, we speak of a {\sl positive g-representation}, or just a {\sl g-representation}.

If $\rho$ is injective, we call $\sigma$ a
{\sl faithful} g-representation.
\end{defi}

\begin{rem} Since this paper deals with $*$-regular rings only, we
  suppress the adjective {\sl positive} when speaking of a linear or a
  generalised representation
  of a $*$-regular ring. We use the term {\sl representation} for a
  linear representation as well as for a generalised representation, if the context leaves no
  ambiguity or both concepts are considered simultaneously.  
\end{rem}

\begin{rem} Note that the properties of a structure to be a {\sl
    (faithful linear) representation of a ring} can be expressed in
    first-order logic \cite{flor}.
\end{rem}

\begin{defi}\lab{drepMOL}{\bf Representation of a (sectional) MOL}\\
A {\sl representation} of a (sectional) MOL $L$ consists of a tuple
$$\varsigma = (D,V_D, \lan \cdot , \cdot \ran, \iota)$$
with $D,V_D, \lan \cdot , \cdot \ran,$ as above and a morphism 
$$\iota: L \ra L(V_D, \lan \cdot , \cdot \ran)$$
of (bounded) lattices between $L$ and the subspace
lattice of $V_D$ such that the (sectional) orthocomplementation
on $L$ corresponds to the (sectional) orthocomplementation on $V_D$
given by the scalar product, that is, for all $x \in L$, we have $\iota(x') =
\iota(x)^\perp$ ($\iota(x^{'_b}) = \iota(x)^\perp \cap \iota(b)$ for a
sectional MOL). 

We call a representation $\varsigma$ {\sl faithful} if the
morphism $\iota$ is injective. $g$-representations are
defined, analoguously.

\end{defi}

A \emph{representation} of an MOL $L$ in an inner product space $(V_F,\Phi)$
is an $(0,1)$-lattice homomorphism $\varepsilon\colon L\to\Sub(V_F,\Phi)$
such that $\varepsilon(x^\prime)=\varepsilon(x)^\perp$ for all $x\in L$. Observe
that, by modularity, $\varepsilon(x)=\varepsilon(x)^{\perp\perp}$ for all $x\in L$.
A representation $\varepsilon$ is \emph{faithful}, if it is one-to-one.
Both  $*$-regular rings and MOLs will be called
{\em representable} if they admit som faithful representation.

\section{The Variety of $*$-Regular Rings}

The term variety is used in the usual sense: A variety is a class of
algebraic structures of the same type that is closed under products,
homomorphic images and substructures. Obviously, an arbitrary product
of $*$-regular rings, where the operations are as usual defined
componentwise, is itself a $*$-regular ring. For homomorphic
images, the following holds. 

\begin{prop}\lab{PHomReg}A homomorphic image of a $*$-regular
  ring is $*$-regular. 
\end{prop}

\begin{bew}Due to \cite{good}, Lemma 1.3 and \cite{flor}, Proposition
  1.7, every two-sided ideal of a
  $*$-regular ring is $*$-regular.
\end{bew}\\

For substructures, we recall the notion of the {\sl Rickart relative inverse}
of an element of a $*$-regular ring. Some preliminary work is needed.

\begin{defi}\lab{dleftrightprojection}{\bf Left and right projection}\\
Let $R$ be a $*$-regular ring. For an element $a \in R$, we call the
unique projection $e$ in $R$ that generates the principal right ideal $aR$ the
{\sl left projection of $a$} and the unique projection $f$ in $R$ that
generates the principal left ideal $Ra$ the {\sl right projection of $a$}.
\end{defi}

\begin{rem}This terminology can be found in \cite{kap},
  p. 27--28 or \cite{KapContGeo}, p. 525. 
We denote the left and right projection of $a$ by $l(a)$ and $r(a)$,
respectively. Furthermore, if $R$ has a unit, we have 
$$ann^l_R(a) =  R(1-e) \quad \mbox{ and } \quad ann^r_R(a) = (1-f)R.$$
\end{rem}

The following result holds.

\begin{lem}\lab{lleftrightprojection} The left and right projection
  of an element $a$ can be constructed in the following way:
For $x \in R$, we set
$$l(x):=x(x^*x)'x^* \quad \mbox{ and } \quad r(x):=x^*(xx^*)'x,$$
where $x'$ denotes any quasi-inverse of $x$.
\end{lem}

\begin{bew2}See \cite{flor}, p. 9--10.
\end{bew2}

\begin{lem}\lab{lrelativeinverse} Let $R$ be a $*$-regular
  ring. Then for each element $a \in R$ there exists a unique element
 $q(a)$ such that the following conditions hold.
\begin{enumerate}
\item $e:= l(a) = a q(a^*a)a^*$.
\item $f:= r(a) = a^*q(aa^*)a$.
\item $fq(a) = q(a)$.
\item $aq(a) = e$.
\end{enumerate}
Furthermore, $q(a)$ has the properties that $q(a)a=f$, the left
projection of $q(a)$ is $f$ and the right projection of $q(a)$ is $e$.
\end{lem}

\begin{bew2} See \cite{kap} or \cite{KapContGeo}). We have defined a
  function $q:R \ra R$ that maps each $a \in R$ to the unique element
  $y$ with the listed properties.
\end{bew2}  

\begin{rem} We call $q(a)$ the {\sl relative inverse of
$a$}. We note that $a$ is the relative inverse of $q(a)$, so
$q^2=id_R$. 
\end{rem}

We arrive at the following result.

\begin{prop}\lab{PSubReg} Let $R$ be a $*$-regular ring.\\ 
The q-subrings of $R$ are exactly the 
$*$-regular subrings of $R$. 
Consequently, we incorporate the unary map $q: R \ra R$
into the signature of $*$-regular rings, that is, a $*$-regular ring $R$
is an algebra of type $(R,+, \cdot, ^*, q, 0)$.   
\end{prop}

\begin{bew} 
Assume that $S$ is closed under $q$. For an element
  $a \in S$, the map $q$ gives a quasi-inverse $q(a)$ of $a$, so $S$
  is regular. Since $S$ is a $*$-subring of the $*$-regular subring
  $R$, $S$ is itself $*$-regular.
\smallskip\\
Conversely, assume that $S \leq R$ is a $*$-regular
  subring of $R$. Let $x \in S$. Due to Lemma \ref{lleftrightprojection},
  we can construct the left and the right projection of $x$ within
  $S$, using the involution on $S$ and any quasi-inverses of $x, x^*,
  xx^*, x^*x$ in $S$. By Lemma \ref{lrelativeinverse}, there exists an
  element $y$ with the desired properties within the $*$-regular ring
  $S$. Since $y$  is the unique element with these properties, we have
  $y=q(x)$. Hence, $S$ is closed under $q$.  
\end{bew}\\

Combining Propositions \ref{PHomReg} and \ref{PSubReg}, we have proven
the first result.

\begin{thm}The class ${\cal R}$ of $*$-regular rings forms a variety. 
\end{thm}

\subsection{Directed Unions and Rings without Unit}

\begin{defi}\lab{ddirectedunionofrings}{\bf Directed union of rings}\\
Let $R$ be a ring and ${\cal S} = \{S_i: i \in I\}$ be a directed family of
  subrings of $R$. We say that $R$ is a {\sl directed union} of the
  family ${\cal S}$ if for each $r \in R$ there exists $k \in I$ such
  that $r \in S_k$.
\end{defi}

\begin{rem}
Casually, we speak of {\sl $R$ being the directed union of the $S_i$}, without giving
the family of the $S_i$ an extra name, and we write $R = \bigcup_{i \in
  I} S_i$, using the usual symbol for an ordinary union. Of course, an arbitrary
union of rings is in general not a ring; hence, the lax notion does
not lead to the risk of misunderstandings.
\end{rem}

\begin{lem}\lab{lultraprodofsubrings} Let $R$ be a $*$-regular ring and assume
  that $R$ is the directed union of a family ${\cal S}$ of $*$-regular
  subrings $S_i$ of $R$. Then $R$ is a $*$-regular subring of an
  ultraproduct of the rings $S_i$, $i \in I$.
\end{lem}

\begin{bew2} Since the class of all $*$-regular rings forms a
  variety, this follows from \cite{gorb}, Theorem 1.2.12 (1).
\end{bew2}

\subsection{Representability and Universal Algebra}

We finish the first section with a look on representability of
$*$-regular rings under an universal-algebraic perspective. 

\begin{lem}\lab{lsubringrep} Let $R$ be a $*$-regular ring with a
  representation $\sigma=(D,V_D,\lan \cdot , \cdot \ran, \rho)$. Then
  every $*$-regular subring $S$ of $R$ is representable. If the
  representation of $R$ is faithful, so is the representation of $S$.  
\end{lem}

\begin{bew2}Just take the restriction $\rho_{|_S}$. 
\end{bew2}

\begin{prop}\lab{pultraprodrep}Let $\{S_i: i \in I\}$ be a family of
$*$-regular rings, $I$ an arbitrary index set. Assume that each $S_i$ 
  admits a linear positive representation. Let $U$ be an ultrafilter
  on $I$. 

Then the ultraproduct  
$$R:= (\prod\limits_{i \in I} S_i)/U$$
admits a linear positive representation. If every $S_i$ has a {\em faithful}
linear positive representation, then so does $R$. 
\end{prop}

\begin{bew} Consider the class of 2-sorted structures 
\begin{eqnarray*} {\cal K}:= \{ (R,V): && R \mbox{ a $*$-regular ring,
  } V  \mbox{ a vector space such that }
\\&& R \mbox{ has a linear positive representation in } V \} 
\end{eqnarray*}

Since the relation that the $*$-regular ring $R$ has a (faithful) linear
positve representation in the vector $V$ can be expressed in
first-order logic, an ultraproduct of a family of structures
$(R_i,V_i) \in {\cal K}$ lies again in ${\cal K}$. 
\end{bew}

\section{Representability of $*$-Regular Rings}

This section is devided into the following parts: In the first part,
we will introduce notation and convention. In the second
part, we will develop the general framework that is needed to tackle
the problem of representability. In the third part, we will present
a proof that a  $*$-regular ring $R$ is g-representable
if and only if so is $\LR$.

\subsection{Convention and Notation}

We consider right modules over rings, denoted by $M_S,
N_T$. Submodules will be denoted by $M_i, N_i$, neglecting the
respective underlying ring. If the contrary is not explicitly stated
(or obvious from the context), we assume that the underlying ring is a
$*$-regular ring (with or without unit). 

For morphisms between submodules $M_i,M_j \leq M$, we write
$\varphi_{ji}: M_i \ra M_j$, where the indices should be read from
right to left. If $M_i,M_j$ have trivial intersection, we define the
graph of a morphism $\varphi_{ji}: M_i \ra M_j$ by
$$\Gamma(\varphi_{ji})= \{x - \varphi_{ji}(x): x \in M_i\}.$$ 

\begin{obs}Let $M_i \cap M_j = \0$. Note that $\Gamma(\varphi_{ji})$
  is a relative complement of $M_j$ in $[0, M_i + M_j]$ and,
  conversely, each such relative complement gives rise to a morphism
  $\psi_{ji}: M_i \ra M_j$. 
\end{obs}

\begin{obs}\lab{naturalsense} Let $M$ be a module with a direct
  decomposition 
$$M = \bigoplus\limits_{i \in I} M_i$$
and denote the corresponding projections and  embeddings by $\pi_i$
and $\vareps_j$, respectively.

Consider a morphism $\varphi_{ji}: M_i \ra M_j$. Then the composition
of $\varphi_{ji}$ with the projection $\pi_i: M 
\ra M_i$ yields a morphism  $\varphi_{ji} \circ \pi_i: M \ra M_j$
defined on all of $M$, i.e.,, $\varphi_{ji} \circ \pi_i \in
Hom(M,M_i)$. Since $M_i \leq M$, we can consider $\varphi_{ji} \circ
\pi_i$ as an element of $End(M)$, too. For the latter point of view, the
formally correct approach would be to consider $\vareps_j \circ
\varphi_{ji} \circ \pi_i$. To avoid technical and notational overload,
we will treat $\varphi_{ji} \circ \pi_i$ itself as an element of $End(M)$.
Note that the composition $\varphi_{ji} \circ \pi_i$ is nothing else
than the extension of the map $\varphi_{ji}:M_i \ra M_j$ to the module
$M$, by defining the action of the extension to be trivial on the other
summands of $M$. Very rarely, we write $\overline{\varphi_{ji}}$ for this
  extension: We just use overlined symbols if we want to distinguish
  between a partial map and its extension. 

Conversely, consider a morphism $\varphi \in M$. We define 
$$\varphi_i:= \varphi \circ \eps_i: M_i \ra M \qquad \varphi_{ji}:=
\pi_j \circ \varphi_i = \pi_j \circ \varphi \circ \eps_i: M_i \ra
M_j$$
Then we have a 1-1-correspondence between a
morphism $\varphi: M \ra M$, and a family $\{\varphi_i: i \in I\}$,
where $\varphi_i: M_i \ra M$, and a family $\{\varphi_{ji}: i,j \in
I\}$, where $\varphi_{ji}: M_i \ra M_j$ since each $\varphi \in
End(M)$ can be decomposed in the following ways:
$$ \varphi = \bigoplus_{i \in I} \varphi_i = \bigoplus_{i \in I}
\sum\limits_{j \in I}\varphi_{ji}$$ 
We agree to write $\varphi = \sum_{i,j \in I} \varphi_{ji}$, with the
convention stated above. We agree to not impose a rigorous notational
strictness, but to understand the notation in the natural sense.  
Similar to the observation above, we note that $\varphi_i
\circ \pi_i$ is nothing else that the extension of the map $\varphi_i =
\varphi \circ \eps_i: M_i \ra M$ to all of $M$.

We note that the conventions are compatible with addition and
multiplication: We can form the sum $\varphi_{ji} +
\psi_{ji}$ and the composition $\varphi_{jk} \circ \psi_{ki}$ in the
natural sense, and for $\varphi, \psi \in End(M)$,
we have 
$$(\varphi + \psi)_{ji} = \varphi_{ji} + \psi_{ji} \quad \mbox{ and }
\quad (\varphi
\circ \psi)_{ji} = \sum_{k \in I} \varphi_{jk} \circ \psi_{ki}.$$

We agree to write $1 = id_M:M \ra M$. Then we have 
$1_{ii} = id_{M_i}: M_i \ra M_i$ (that is, the corresponding
extension $1_{ii}$ acts like the identity on $M_i$ and trivially on
every other summand $M_j$) and $1_{ji} = 0_{ji}: M_i \ra M_j$ (that is,
the extension $1_{ji}$ coincides with the zero map on $M$).
\end{obs}

\begin{obs}For cyclic modules $M_S,N_S$ with
  generators $x,y$, a morphism $\varphi: M \ra N$ is determined by its
  action on the generator $x$ of $M$. If $M=xS$, we have $f(xs) = 
  f(x)s$ for every $xs \in M$. Conversely, each choice $z \in yS$
  defines a morphism $g: M \ra N$ via $xs \mapsto zs$. 

In particular, let $R$ be a regular ring and consider the module
  $R_R$. Assume that $I = eR, J = fR$ are principal right
  ideals in $R$ (that is, cyclic submodules of $R_R$). 
Since $R$ is regular, the generators $e,f$ can be taken to be
  idempotent. 

Let $r \in R$ such that $re \in J$, that is, $re = fc$ for some
  $c \in R$ (or, equivalently, $f(re) = re$). Then the left
  multiplication with $r$ defines a right-$R$-module-homomorphism
  $\widehat{r}$ between $I$ and $J$
$$\widehat r: I \ra J \qquad es \mapsto r(es).$$
\end{obs}

\begin{rem}From now on, if possible, we denote the action defined by
  left multiplication with an element $r$ by $\widehat r$. We will speak
  of the {\sl left multiplication morphism} (or {\sl left
  multiplication map} or {\sl  left multiplication}) $\widehat{r}$.
\end{rem}

\subsection{General Framework}

In this section, we will develop the necessary machinery for the proof
of the desired result. In order to simplify the lines of argument and to
clarify the applied technique, we have chosen to separate
the ring-theoretical aspects, the lattice- and frame-theoretical
aspects and the general module-theoretical mechanisms as far as
possible.

\subsubsection{Decomposition Systems \& Abstract Matrix Rings}

\begin{defi}\lab{dDS}{\bf Decomposition system of a module}\\
Let $M_S$ be a right $S$-module over $S$ and $I=\{i: 0 \leq i < n+k\}$
an index set, where $n < \omega$ and $k \leq \omega$. A {\sl decomposition system $\vareps$ of $M$ of format $(n,k)$} consists of 
\begin{enumerate}
\item a decomposition $M=\bigoplus_{i \in I} M_i$ of $M$ into a direct sum
  of submodules, 
\item corresponding projections $\pi_i: M \thra M_i$ and
  embeddings $\eps_i: M_i \hra M$,
\item a family $\{\eps_{ij}: i,j \in I\}$ of maps $\eps_{ij}$,
\item submodules $z_{ij}$ of $M$ for $i \in I, j < n$, and 
\item a 1-subring $C \leq End(M_0)$
\end{enumerate}
such that the following conditions are satisfied:
\begin{enumerate}
\item For $i=j$, we have $\eps_{ii} = id_{M_i}$. 
\item For $i,j < n$, $\eps_{ij}, \eps_{ji}$ are mutually inverse
  morphisms, i.e., $\eps_{ij} \circ \eps_{ji} = id_{M_i}$. 
\item For $i \in I$, we have $\eps_{i0} \circ \eps_{0i} =
  id_{M_i}$ (in particular, $\eps_{0i}$ is injective).
\item For distinct indices $i,j,k$ such that $k,j < n$, we have
  $\eps_{ki} = \eps_{kj} \circ \eps_{ji}$
\item For $j < n$, $z_{ij}$ is a relative complement of
  $im(\eps_{ji})$ in $[0,M_j]$.
\item For $i \in I$, $\eps_{0i} \circ \eps_{i0} \in C$.
\end{enumerate}

In other words, for $i,j < n$, the submodules $M_i,M_j$ are
isomorphic, while for $i \in I, j < n$, $M_i$ is isomorphic to a
submodule of $M_j$ -- and the morphisms $\eps_{ji}$ are the
corresponding isomorphisms and embeddings. 

The relative complements
$z_{ij}$ are integrated into the notion of a decomposition systems
for the following reason: For $i,j$ with $j < n$, the injective
morphism $\eps_{ji}:M_i \hra M_j$ has a left inverse
$\eps_{ij}:: M_j \ra M_i$, defined only on $im(\eps_{ji}) \leq
M_j$. Taking the relative complement $z_{ij} \leq M_j$ of
$im(\eps_{ji})$ in $[0,M_j]$, we can extend the partial morphism
$\eps_{ij}:: M_j \ra M_i$ to a morphism $\eps_{ij}: M_j \ra M_i$
by setting $\eps_{ij}(x):= 0$ for all $x \in z_{ij}$ (i.e., the
extension $\eps_{ij}: M_j \ra M_i$ acts trivially on $z_{ij}$).
\end{defi}

\begin{rem} For the ease of notation, we stated that a decomposition
  system contains a family of maps $\eps_{ij}$ for $i,j \in
  I$. The required conditions should have made clear that
  only particular maps have to exist. Of course, the maps that do
  exist are (partial) morphisms satisfying the desired relations. 
(One might take the view that the {\sl other} maps are partial maps
  with trivial domain.)

We write $\vareps = \vareps(C,M)$ to indicate the ring
  $C$ and the module $M$ under consideration.

We recall Observation \ref{naturalsense} for the natural
identifications and conventions.
\end{rem}

\begin{defi}\lab{dMorphDSs}{\bf Morphisms between decomposition
    systems}\\
Let $M_S,M_{S'}'$ be modules over $S,S'$ and $\vareps, \vareps'$
decomposition systems of $M,M'$, respectively.\footnote{Similarly, we
  denote the components of the two systems by the same letters, once
  with prime, once without.}
A morphism between the two decomposition systems $\vareps, \vareps'$
is a map $\eta:\vareps \ra \vareps'$ such that the components of
$\vareps$ get mapped onto the components of $\vareps'$. In particular,
the following hold.
\begin{enumerate}
\item $\eta(M_i) = M_i'$ and $\eta(\pi_i) = \pi'_i$, $\eta(\eps_i) =
  \eps'_i$ for all $i \in I$. 
\item $\eta(z_{ij}) = z'_{ij}$ for all $i,j \in I$.
\item $\eta(\eps_{ij}) = \eps'_{ij}$ for all $i,j \in I$.
\item $\eta: C \ra C'$ is a morphism of rings with units.
\end{enumerate}

A morphism $\eta$ between decomposition systems will be called {\sl
  injective} or an {\sl embedding of decomposition systems} if $\eta: C
  \ra C'$ is injective.  
\end{defi}

\begin{defi}\lab{dabstractmatrixring}{\bf Abstract matrix ring}\\
Let $M_S$ be a module and $\vareps$ a decomposition system
  of $M$. The {\sl abstract matrix ring} with respect to the
decomposition system $\vareps$ of $M$ is
$$R(\vareps, C, M):=\{\varphi \in End(M_S): \eps_{0j} \circ \varphi_{ji}
\circ \eps_{i0} \in C \mbox{ for all } i,j \} \sub End(M_S)$$
where, as above, $\varphi_{ji} = \pi_j \circ \varphi \circ \eps_i$ and
$\pi_j, \eps_i$ are the natural projections and embeddings
belonging to decomposition system $\vareps$.
\end{defi}

The following result justifies this definition. 

\begin{prop}\lab{pabstractmatrixring} The set $R(\vareps, C,M)$ is a
  1-subring of $End(M_S)$.
\end{prop}

\begin{prop}\lab{pextensioneta}Let $M_S,M'_{S'}$ be two modules with
  decomposition systems $\vareps, \vareps'$ and $\eta: \vareps \ra
  \vareps'$ a morphism of decomposition systems between $\vareps$ and
  $\vareps'$. 
Declaring
$$\eta(\varphi_{ji}):= \eps'_{j0} \circ \eta\Big( \eps_{0j} \circ
\varphi_{ji} \circ \eps_{i0} \Big) \circ \eps'_{0i}$$ for morphisms
$\varphi_{ji}: M_i \ra M_j$, the map $\eta$ can be extended to a map 
$$\eta: R(\Phi,C,M) \ra R(\Phi',C',M')$$
in the following way. Since $\varphi \in End(M)$ decomposes into $\varphi = \bigoplus \varphi_i =
\sum \varphi_{ji}$, we 
can define $ \eta(\varphi_i):= \sum_{j \in I}  \eta(\varphi_{ji})$ for
a fixed index $i \in I$ and 
$$\eta(\varphi):= \sum_{i \in I} \eta(\varphi_i) = \sum_{i,j \in I}
\eta(\varphi_{ji}).$$ 
With this definition, $\eta: R(\Phi,C,M) \ra R(\Phi',C',M')$ is a
morphism of rings with unit. If the restriction $\eta_{|_C}: C \ra C'$
is injective, then so is the map $\eta:  R(\Phi,C,M) \ra R(\Phi',C',M')$.
\end{prop}

\subsubsection{Frames and Induced Structures}

In this section, we approach the connection between the general
framework and our particular setting. Starting with a frame in $L(M)$,
we will develop the notion of the coefficient ring of a frame and the
notion of an induced decomposition system.

\begin{defi}\lab{dcoefficientring}{\bf Coefficient ring of a
    frame}\\
Let $M_S$ be a right module over
  $S$ and $\Phi$ a stable $(n,k)$-frame in $L(M_S)$, contained in the
  sublattice $L \leq L(M_S)$ with $n \geq 3$. The {\sl coefficient
  ring} of $(\Phi,L,M)$ is 
$$C(\Phi, L, M):=\{\varphi \in End(M_0): \Gamma(\eps_{10} \circ \varphi)
  \in L\} \sub End(M_0)$$
\end{defi}

The following result justifies this definition.

\begin{prop}\lab{pcoeffcientring} The set $C(\Phi, L, M)$ is a
  1-subring of $End(M_0)$.
\end{prop}

\begin{rem}The lines of argument and the technique of this proof are
  well-known: It is possible to express the ring operations via lattice
  terms with constants in $\Phi$. (See the works of von Neumann,
  J\'onsson and Handelman.) These terms are
  uniform in the frame $\Phi$. In particular, they are independent of the
  particular module $M_S$.
\end{rem}

\begin{prop}\lab{pDSofaframe}{\bf The Decomposition System of a
  Frame}\\
 Let $M_S$ be a right module over $S$ and
  $\Phi$ a stable $(n,k)$-frame in $L(M_S)$ contained in the
  sublattice $L \leq L(M_S)$ with $n \geq 3$.

Then $\Phi$ induces a
  decomposition system in the following way.
Since $\Phi$ is a frame, we have a decomposition of $M$ into a
direct sum $M=\bigoplus M_i$, together with corresponding projections
and embeddings. As usual, the axes of perspectivity as well as the
axes of subperspectivity are the graphs of
morphisms between the summands. Since $\Phi$ is stable, we have
relative complements $z_{ij}$ as required.\footnote{That is,  $z_{ij}$
  a relative complement of $b_{ji}$ in $[0,a_j]$, where $b_{ji}$ is
  the image of $a_i$ under  $a_{ji}$ in $[0,a_j]$.}  As ring $C$, we
take the coefficient ring $C(\Phi,L,M)$.\smallskip   

We denote the decomposition system induced by $\Phi$ 
by $\xi = \xi_{\Phi,L}(C,M)$. 
\end{prop}

\begin{bew} The only thing left to show is Property 5 in
  Definition \ref{dDS}. Consider $\eps_{0i}, \eps_{i0}$. Both graphs
  $\Gamma(\eps_{0i}), \Gamma(\eps_{i0})$ and of course
  $\Gamma(\eps_{10})$ are part of the frame $\Phi$. Since we can
  express composition of maps  by lattice terms 
  with constants in $\Phi$, we have $\Gamma(\eps_{10} \circ \eps_{0i}
  \circ \eps_{i0}) \in L$. 
\end{bew}

\begin{defi}\lab{dmatrixringofaframe}{\bf Matrix ring of a frame}\\
Let $M_S$ be a right module over $S$ and $\Phi$ a stable
  $(n,k)$-frame in $L(M_S)$ (with $n \geq 3$), contained in the
  sublattice $L \leq L(M_S)$, $C(\Phi, L, M)$ the coefficient ring as
  defined in Defintion \ref{dcoefficientring} and $\xi=\xi_{\Phi,L}(C,M)$ the
  induced decomposition system as defined in Definition
  \ref{pDSofaframe}. The ring $$ R(\Phi,L,M):= R(\xi,C(\Phi,L,M),M)$$
will be called the {\sl matrix ring (of $\Phi,L,M$)}. 
\end{defi}

We consider the following situation: Let $M$ and $M'$ be modules over
$S$ and $S'$, $L \leq L(M_S)$ a complemented $0$-$1$-sublattice and
$\Phi$ a stable frame in $L(M_ S)$
contained in $L$ of format $(n,k)$ with $n \geq 3$. Assume that we are given a
morphism $\iota: L  \hra L(M')$ of bounded complemented lattices.

\begin{obs}\lab{oimageframe} The image $\Phi':=\iota[\Phi]$ is a
  stable frame in $L(M')$, contained in $L':=\iota[L] \leq L(M')$. In particular, 
  we have $\iota(M_i)= M'_i$, $\iota(\pi_i)=
  \pi'_i$, $\iota(\eps_i)= \eps'_i$ and $\iota(z_{ij}) = z_{ij}$.
\end{obs}

\begin{prop}\lab{pinducedeta} The morphism $\iota: L \ra L'$ induces
  a morphism $\eta$ between the induced decomposition systems
  $$\xi:=\xi_{\Phi,L}(C(\Phi,L,M),L,M) \quad \mbox{ and } 
\quad   \xi':=\xi'_{\Phi',L'}(C(\Phi',L',M),L',M').$$

If $\iota: L \ra L'$ is injective, then so is $\eta: \xi \ra \xi'$.  
\end{prop}

\begin{bew} We want to define $\eta$ via the lattice morphism $\iota:
  L \ra L'$.
For the first two properties of a morphism between two decomposition
systems (see Definition \ref{dMorphDSs}), we define $\eta$ to coincide
with $\iota$ on   the submodules $M_i, z_{ij}$ of $M$ and recall Observation
  \ref{oimageframe}. 

Now consider the morphisms
$\eps_{ji}$ given by the frame $\Phi$, i.e., $\Gamma(\eps_{ji}) =
a_{ji} \in \Phi$. Then  
$$\iota\big(\Gamma(\eps_{ji})\big) = \iota(a_{ji})  = a'_{ji} =
\Gamma(\eps'_{ji}) \in \Phi'$$
Setting
$$\eta(\eps_{ji}):=\eps'_{ji} \mbox{ for }i \neq j < n  \quad
\mbox{and} \quad 
\eta(\eps_{ii}) : = \eps'_{ii} \mbox{ for  arbitrary } i $$ 
we have guaranteed that $\eta$ maps the morphism $\eps_{ji}$ to
$\eps'_{ji}$.

For appropriate indices $i,j,k$, the compatibility $\eps_{ki} = \eps_{kj}
\circ \eps_{ji}$ is determined by the lattice-theoretical equation
$$[a_{kj} + a_{ji}] \cdot [a_k + a_i] = a_{ki}$$
of the elements of the frame $\Phi$ (and similarly for $\Phi'$). Therefore, 
we have 
$$ \eta(\eps_{kj} \circ \eps_{ji}) = \eta(\eps_{ki}) =
\eps'_{ki} = \eps'_{kj} \circ \eps'_{ji}$$
for   appropriate indices $i,j,k$.

Secondly, we consider an element $\varphi$ of the coefficient ring $C
= C(\Phi,L,M)$, i.e., $\varphi: M_0 \ra M_0$ with
$\Gamma(\eps_{10} \circ \varphi) \in L$.

The property that $\eps_{10} \circ \varphi$ is a morphism between
$M_0$ and $M_1$ is equivalent to the lattice-theoretical property that 
$\Gamma(\eps_{10} \circ \varphi)$ is a relative complement of $M_1$
in $[0,M_0 + M_1]$. Since $\iota: L \ra L'$ is a
lattice morphism mapping $\Phi$ to $\Phi'$, $\iota(\Gamma(\eps_{10}
\circ \varphi))$ is a relative complement of $M'_1$ in $[0,M'_0 + M'_1]$,
i.e., the graph of a morphism $\psi: M'_0 \ra M'_1$. Composing $\psi$
with $\eps'_{01}$, we can define $$\eta(\varphi):= \eps'_{01} \circ
\psi: M'_0 \ra M'_0$$

Thirdly, we can capture the ring operations on $C(\Phi,L,M)$ via
lattice terms with constants in $\Phi$. Hence, the ring operations are
transferred via $\iota: L \ra L'$ to $\Phi'$ and $C'$. Accordingly, the map 
$\eta: C(\Phi,L,M)  \ra C(\Phi',L',M')$ is  a morphism of rings. 

Finally, injectivity of $\iota$ implies injectivity of $\eta$.
\end{bew}

\begin{cor}\lab{cringmorphismeta}
In the given situation, there exists a morphism 
$$\eta: R(\Phi,L,C) \ra R(\iota[\Phi], \iota[L], C')$$
of rings with unit.

If $\iota: L \ra L'$ is injective, then so is $\eta: R \ra R'$. In
particular, if $L$ embedds into the subspace lattice $L(V)$ of a 
vector space $V$, we have a ring embedding
$$\eta: R(\Phi,L,C^M) \hra End(V_D)$$
\end{cor}

\begin{bew2} Combine Proposition \ref{pinducedeta} and  Proposition
  \ref{pextensioneta}.
\end{bew2}

\subsection{Representability of $*$-Regular Rings}

This section is dedicated to the desired result on representability of
$*$-regular rings $R$ such that $\LR$ is representable. First,
we focus our attention on a $*$-regular ring $R$ with unit such that 
the MOL $\LR$ contains a stable orthogonal frame. With that
restriction, we aim at representability of simple $*$-regular rings
with unit. Subsequently, we will deal with simple $*$-regular rings
without unit and finally, with subdirectly irreducible $*$-regular
rings (with and without unit). Due to the first main theorem that the
class of all $*$-regular rings is a variety, with Theorem
\ref{csi*regringrep} , we reach the desired result that a
$*$-regular ring $R$ is g-representable 
if so is  $\LR$.

\subsubsection{$*$-Regular Rings with Frames}\lab{ssRegRingsFrames}

\begin{rem}\lab{rAssumptions} For this subsection, we assume that 
\begin{enumerate}
\item $R$ is a $*$-regular ring with unit,
\item $L:= \LR$ is a MOL of height $h(L) \geq 3$,
\item $\Phi$  is a stable orthogonal $(n,k)$-frame in $L$ with $n
  = 3$, 
\item  $M = R_R$, if not stated otherwise.
 \end{enumerate}
\end{rem}

\begin{rem}\lab{rConsequenceAssumptions}  Moreover, we assume that 
there exists a faithful representation $\iota: L \hra
  L(V_D, \lan \cdot , \cdot \ran)$ of the MOL $L$. Consequently,
  Corollary \ref{cringmorphismeta} applies in its full strength. 
\end{rem}

\begin{cor}\lab{ccyclic}Let $e_i, e_j$ be projections in $R$ and $e_iR, e_jR$
  the corresponding cyclic modules. Any morphism $\varphi_{ji}: e_i R
  \ra e_jR$ is a left multiplication by a ring element $e_j s
  e_i \in e_j R e_i$. 
\end{cor}

\begin{obs}\lab{oinducedringDS} By Proposition \ref{pDSofaframe},
  the stable orthogonal frame $\Phi$ induces a decomposition system $\xi =
  \xi_{\Phi,L}(C,M) =  \xi_{\Phi,L}(C,R_R)$. More exactly, we have the
  following correspondences. 
\begin{enumerate} 
\item The summands $M_i$ correspond to principal right ideals $e_iR$
  generated by a projection $e_i$. Each projection $\pi_i$ corresponds
  to a map $\widehat{e_i}$ given by left multiplication with $e_i$
  and coincides with the (extension of the) embedding $\eps_i =
  id_{M_i}$.  
\item The morphisms $\eps_{ji}: e_iR \ra e_jR$ are given by
  $\widehat e_{ji}$ with $e_{ji}:=\eps_{ji}(e_i)$ an element of $e_i R e_j$.
\item For the coefficient ring of the frame, we have 
$$C = C(\Phi,L,M) = C(\Phi,L,R_R) = \{\widehat{r}: r \in e_0 R e_0\}.$$
\end{enumerate}
\end{obs}

\begin{rem} From now on, we will denote the morphisms $\eps_{ji}$ given
  by the decomposition above by $\eps_{ji}$ and $\widehat{e_{ji}}$
  interchangedly, as it suits the particular situation.
\end{rem}

\begin{cor}\lab{cringisomorphismtheta}We have an isomorphism 
$$\theta: R \ra R(\Phi,L,R_R)$$
of rings with unit.
\end{cor}

\begin{bew2} As stated in Remark \ref{rAssumptions} at the beginning of the
  section, we agree to write $M$ for $R_R$. We recall the definition
  of the matrix ring of a frame:  
$$ R(\Phi,L,M) = \{ \varphi \in End(M): \forall i,j \in I. \;
\eps_{0j} \circ \varphi_{ji} \circ \eps_{i0} \in C(\Phi,L,M)\}$$
Since $R$ contains a unit, an endomorphism $\varphi \in End(M)$ is
given by left multiplication $\widehat{r}$ for some $r \in R$. We notice that 
$$(\widehat{r})_{ji} = \pi_j \circ \widehat{r} \circ \eps_i =
\widehat{e_j} \circ \widehat{r} \circ \widehat{e_i} = \widehat{e_j r
  e_i} = \widehat{r_{ji}}, 
\quad \mbox{ with } r_{ji}:=  e_j r e_i.$$
Then we have
$$\eps_{0j} \circ (\widehat{r})_{ji} \circ \eps_{i0} = \hat e_{0j} \circ
\widehat{e_j r e_i} \circ \hat e_{i0} = \widehat{e_{0j} r  e_{i0}},
\quad \mbox{ since } e_{0j} e_j = e_{0j}, e_i e_{i0} =   e_{i0}.$$
The equality 
$$\eps_{0j} \circ (\widehat{r})_{ji} \circ \eps_{i0} =
\widehat{e_{0j} r  e_{i0}}$$
and Observation \ref{oinducedringDS} lead
to $\eps_{0j} \circ (\widehat{r})_{ji} \circ \eps_{i0} \in
C(\Phi,L,M)$. Since the indices $i,j$ were arbitrary, we have
$\widehat{r} \in R(\Phi,L,M)$. 

In particular, for an element $r \in R$, the left multiplication 
$\widehat{r}: M \ra M$ decomposes into  
$$\widehat{r} = \sum \widehat{r_{ji}} \quad \mbox{ where  }
\widehat{r_{ji}}: e_iR \ra e_jR, \quad \mbox{ and  } r_{ji} = e_jre_i
\in e_j R e_i.$$
that is, the isomorphism $\theta: R \ra R(\Phi,L,M)$ is given by
$\theta: r \mapsto \widehat{r}$.

Of course, we have
\begin{eqnarray*}
\Gamma(\eps_{10} \circ \widehat{e_{0j} r  e_{i0}}) &=&
\Gamma(\widehat{e_{10}} \circ \widehat{e_{0j} r e_{i0}} ) =\Gamma(
\widehat{e_{1j} r e_{i0}})\\ & = & (e_0 - e_{1j}re_{i0})R \, \in L=\LR
\end{eqnarray*}
\end{bew2}

\begin{cor}\lab{cringembeddingrho}We have an embedding 
$$\rho: R \hra End(V_D)$$
of rings with unit.
\end{cor}

\begin{bew2} By Remark \ref{rConsequenceAssumptions} 
  the MOL $L=\LR$ is assumed to be   representable in $L(V_D, \lan .,. \ran)$.
We combine that with Corollaries \ref{cringmorphismeta} and
\ref{cringisomorphismtheta} and define $\rho:=\eta \circ \theta$ to
get the desired isomorphism. 
\end{bew2}\\

It is left to show that the isomorphism translates the involution on 
$R$ into adjunction with respect to the scalar product. In the
following, assume that in addition to the assumptions of ref{...}, we
have the following.

\begin{enumerate}
\item $(V_D, \lan .,. \ran)$ a vector with scalar product,
\item $K$ a MOL represented in $L(V_D, \lan .,. \ran)$, i.e., $K$
  is a modular sublattice of $L(V_D)$ and the orthocomplementation on
  $K$ is  induced by the scalar product $\lan .,. \ran$ on $V_D$,
\item $\Psi$ a stable orthogonal frame in the MOL $K$, 
\item $U_i,U_j \leq V$ with $U_i,U_j \in \Psi$ and $f:U_i \ra U_j, g:
  U_j \ra U_i$ linear maps, 
\end{enumerate}
If we discuss both situations - that is, for the frame $K$ with its
frame $\Psi$ or subspaces of $V$ - simultaneously, we use the symbol
$\Upsilon$ for the frame, $N_i$ for submodules or subspaces and $a,b$
for morphisms. 

\begin{defi}\lab{dadjointnessEndVD}{\bf Adjointness on $End(V_D, \lan
    .,. \ran)$}\\
We call {\sl $f$ and $g$ adjoint to each other} (with 
respect to $\lan .,. \ran$) if
$$\forall v \in U_i, w \in U_j. \qquad \lan fv,w \ran = \lan v, gw\ran.$$
\end{defi}

\begin{rem}Notice that $U_i, U_j$ are elements of the orthogonal frame
  $\Psi$, in particular, if $i \neq j$, then $U_i$ and $U_j$ are
  orthogonal to each other.  
\end{rem}

\begin{lem}\lab{ladjointness}The following conditions are equivalent: 
\begin{enumerate}
\item $f$ and $g$ are adjoint to each other in the sense of Definition
  \ref{dadjointnessEndVD}.
\item The extensions $\overline{f}, \overline{g}: V \ra V$ are adjoint to each
other in the usual sense.
\end{enumerate}
If $i \neq j$, both these conditions are equivalent to $ \Gamma(f)
\perp \Gamma(-g)$.\footnote{Notice that $\Gamma(f), \Gamma(g),
  \Gamma(-g)$ are contained in the MOL $K$.}
\end{lem}

\begin{bew2}The equivalence of the first two conditions is
  immediate. Now, if $i \neq j$, we have  
$$\Gamma(f) = \{v - fv: v \in U_i\} \qquad \Gamma(-g) = \{w + gw: w
\in U_j\}$$
and 
\begin{eqnarray*}&&\Gamma(f) \perp \Gamma(-g) \gdw \lan v - fv, w + gw
  \ran = 0 \mbox{ for all } v \in U_i, w \in U_j \\
& \gdw &\lan v , w \ran + \lan v, gw \ran - \lan fv, w \ran  - \lan fv,
gw \ran = 0 \mbox{ for all } v \in U_i, w \in U_j\\
 & \gdw & \lan v, gw \ran - \lan fv, w \ran  = 0 \mbox{ for all } v \in
 U_i, w \in U_j\\
& \gdw & \lan v, gw \ran  = \lan fv, w \ran  \mbox{ for all } v \in
 U_i, w \in U_j\\
& \gdw & \mbox{$f$ and $g$ are adjoint to each other in the sense of
  Definition \ref{dadjointnessEndVD}},
\end{eqnarray*}
where the terms $\lan v , w \ran, \lan fv,gw \ran$ vanish since
$U_i, U_j$ are orthogonal to each other. \end{bew2}

Now, we derive a result similar to Lemma \ref{ladjointness} for a
$*$-regular ring $R$ and the relation between the involution on $R$
and the orthogonality on $L = \LR$. 
 
\begin{lem}\lab{linvolutiongraphs}The involution on $R$ can be
  captured via the orthogonality on $L$, more exactly, for $a_{ij} \in
  e_i R e_j, b_{ji} \in e_j R e_i$ with $i\neq j$, the following
  conditions are equivalent: 
\begin{enumerate}
\item $a_{ij} = b_{ji}^*$
\item $\Gamma(\widehat{a_{ij}}) \perp \Gamma(- \widehat{b_{ji}})$
\end{enumerate}
\end{lem}

\begin{bew2}
Since $\widehat{a_{ij}}: e_jR \ra e_iR$ and $-\widehat{b_{ji}}: e_iR \ra
e_j$, we have 
$$\Gamma(\widehat{a_{ij}}) = (e_j - a_{ij}e_j)R \qquad 
\Gamma(- \widehat{b_{ji}}) = (e_i + b_{ji}e_i)R$$ 
The orthogonality on $L$ is given by $pR \perp qR :\gdw q^*p =
0$. Calculating yields 
\begin{eqnarray*}
(e_i + b_{ji}e_i)^*\cdot (e_j - a_{ij}e_j) &=&  (e_i + e_i
b^*_{ji})(e_j - a_{ij}e_j)\\
& =&  e_i e_j - e_i a_{ij}e_j + e_i b^*_{ji} e_j - e_i b^*_{ji}
a_{ij}e_j\\
&=& e_i (-a_{ij}+b^*_{ji})e_j = -a_{ij}+b^*_{ji},
\end{eqnarray*}
since $a_{ij}, b^*_{ji} \in e_i R e_j$. Hence
$a_{ij} = b^*_{ji}$ iff $\Gamma(\widehat{a_{ij}}) \perp \Gamma(-
\widehat{b_{ji}})$.
\end{bew2}

\begin{cor}\lab{cuniquenessadjoint}{\bf Uniqueness}\\
Let $(i,j)$ be an arbitrary pair of indices.

A linear map $f: U_i \ra U_j$ has at most one adjoint $g:
  U_j \ra U_i$. Due to this uniqueness, it is legitimate to write $f^*
  = g$ if $f,g$ are adjoint to each other.

Likewise, a map $\widehat{a_{ij}}: e_j R \ra e_i R$ gives rise to a
map $\widehat{b_{ji}}: e_i R \ra e_j R$, namely $\widehat{b_{ji}} = 
\widehat{a_{ij}^*}$. If $i \neq j$, we have $\widehat{b_{ji}} = 
\widehat{a_{ij}^*}$ iff $\Gamma(\widehat{a_{ij}})
\perp \Gamma(- \widehat{b_{ji}})$. 
\end{cor}

\begin{lem}\lab{lexistenceadjointeps} For each morphism $\eps_{ki}: N_i
  \ra N_k$, there exists an adjoint $\eps_{ki}^*: N_k \ra N_i$ in
  $\Upsilon$.
\end{lem}

\begin{cor}\lab{cadjointepskiepsikid} For $i,k \in I$ with $k < n$ we have
  $\eps_{ki}^* \circ \eps_{ik}^* = id_{U_i}$.
\end{cor}

\begin{bew2}We have 
$$\eps_{ki}^* \circ \eps_{ik}^* = (\eps_{ik} \circ  \eps_{ki})^* =
(id_{N_i})^*  = id_{N_i}.$$
\end{bew2}

\begin{rem}Obviously, Lemma \ref{lexistenceadjointeps} and
  Corollary \ref{cadjointepskiepsikid} hold for arbitrary  indices
  $i,k \in I$: Recall that for if $i=k$, we have   $\eps_{ik} =
  \eps_{ii} = \eps_i = id_{N_i}$, which is an hermitian idempotent map.  
\end{rem}

\begin{cor}\lab{cadjointdrag} Let $a: N_i \ra N_j$ and $b: N_j \ra
  N_i$ be as before.  

Then $a$ and $b$ are adjoint  to each other iff
$\eps_{i0}^* \circ b \circ \eps_{1j}^* $ and $\eps_{1j}  \circ a \circ
\eps_{i0}$ are adjoint to each other.  
\end{cor}

\begin{bew2}Assume that $\eps_{i0}^* \circ b \circ \eps_{1j}^* $ and
  $\eps_{1j}  \circ a \circ \eps_{i0}$ are adjoint to each other,
  where adjoint is either understood in the sense of
  Definition \ref{dadjointnessEndVD} or in the sense of Lemma \ref{linvolutiongraphs}. Then 
\begin{eqnarray*}
\eps_{0i}^* \circ (\eps_{1j}  \circ a   \circ \eps_{i0})^* \circ
\eps_{j1}^* =  (\eps_{j1}\eps_{1j} a \eps_{i0}\eps_{0i})^* = a^*
\end{eqnarray*}
and 
\begin{eqnarray*}
\eps_{0i}^* \circ (\eps_{1j}  \circ a   \circ \eps_{i0})^* \circ
\eps_{j1}^* = \eps_{0i}^*  \circ (\eps_{i0}^*  \circ b  \circ
\eps_{1j}^*)  \circ \eps_{j1}^* = (\eps_{0i}^* \eps_{i0}^*) b
(\eps_{1j}^* \eps_{j1}^*) = b, 
\end{eqnarray*}
so $a^*=b$.

Now, assume that $a,b$ are adjoint to each other. Then 
$$(\eps_{1j}  \circ a   \circ \eps_{i0})^* = \eps_{i0}^* \circ a^*
\circ \eps_{1j}^* = \eps_{i0}^* \circ b \circ \eps_{1j}^*.$$
\end{bew2}

\begin{rem}Corollary \ref{cadjointdrag} holds for arbitrary indices
  $i,j$, too. In particular, we can complete Lemma \ref{linvolutiongraphs} and
  Corollary \ref{cuniquenessadjoint} by noting the following:
If $a_{ii}, b_{ii} \in e_iRe_i$, we have
$$a^*_{ii}  = b_{ii} \gdw \Gamma(\eps_{i0}^* \circ \widehat{a_{ii}} \circ 
\eps^*_{1j}) \perp \Gamma(- \eps_{1j} \circ \widehat{b_{ii}} \circ
\eps_{i0})$$ 
\end{rem}

\begin{prop}\lab{p*ringembeddingrho} The map 
$$\rho = \eta \circ \theta: R \ra End(V_, \lan .,. \ran)$$ 
defined in Corollary \ref{cringembeddingrho} is a $*$-ring-embedding
of involutive rings with unit.

\end{prop}

\begin{bew} 
First, we recall that the map $\eta: R(\Phi,L,M) \hra End(V_D)$
defined in Corollary \ref{cringmorphismeta} was defined via
the lattice embedding $\iota: L  \hra L(V_D)$. In this situation, we
consider a MOL $L$ and a MOL-embedding of $L$ into $L(V_D, \lan
.,. \ran)$. As shown before, for morphisms $\eps_{ji}:M_i \ra M_j \in
\Phi$, we can define an adjoint 
operator $\eps_{ji}^*: M_j \ra M_i$ via the orthogonality on $L$. For
morphisms $\varphi_{10}: M_0 \ra M_1$ and $\psi_{01}:M_1 \ra M_0$, the
relation of adjointness could be captured via orthogonality on graphs.
In particular, we have $\rho(\eps_{ji}^*)= \big(\rho(\eps_{ji})\big)^*$
and $\rho(\varphi_{10}^*) = \big(\rho(\varphi_{10})\big)^*$.

Now, for $r \in R$, consider $e_j r e_i \in e_j R e_i$.
Then 
\begin{eqnarray*}
 \big(\rho(e_j r e_i) \big)^* &=& \big(\rho(e_{j1} ( e_{1j} e_j r e_ i
 e_{i0}) e_{0i})\big)^* = \big(\rho(e_{j1}) \rho( e_{1j} e_j r e_ i
 e_{i0}) \rho(e_{0i})\big)^* \\
& = & \big(\rho(e_{0i})\big)^* \big(\rho( e_{1j} e_j r e_ i
 e_{i0})\big)^* \big(\rho(e_{j1})\big)^*\\
& = & \rho(e_{0i}^*) \rho\big( (e_{1j} e_j r e_ i
 e_{i0})^*\big) \rho(\eps_{j1}^*)\\
& = & \rho(e_{0i}^*) \rho\big( e_{i0}^* e_i^* r^* e_ j^*
 e_{1j}^*\big) \rho(\eps_{j1}^*)\\
& = & \rho\big(e_{0i}^* e_{i0}^* e_i^* r^* e_ j^*
 e_{1j}^*\eps_{j1}^*\big) = \rho(e_i^* r^* e_ j^*) = \rho\big( (e_j r
 e_i)^* \big).
\end{eqnarray*}

Hence, we have
\begin{eqnarray*}\big(\rho(r)\big)^* & = & \Big(\rho \big(\sum e_j r
  e_i\big)\Big)^* = \Big( \sum \rho(e_j r e_i) \Big)^* = \sum
  \big(\rho(e_j r e_i)\big)^*\\  
& = & \sum \rho\big((e_j r e_i)^*\big) = \sum \rho( e_i^* r^* e_j^*) =
  \rho \Big(\sum e_i^* r^*  e_j^*\Big)\\ & = & \rho\Big(\sum
  (e_jre_i)^* \Big)  =  \rho\Big( \big(\sum e_j r e_i\big)^*\Big)
  = \rho(r^*). 
\end{eqnarray*}
\end{bew}

\subsubsection{Simple $*$-Regular Rings}

\begin{cor}\lab{csimple*regringwithunitisrep} Every simple
  $*$-regular ring $S$ with unit admits a faithful linear
  representation provided that $\LS$ does so.
 \end{cor}

\begin{bew2} We may assume that $S$ is non-Artinian, hence, we can
  assume that $S$ has height at least 
  3. By Corollary \ref{cRsimpleimpliesLRcontainssonkframe}, the
  MOL $L=\LS$ contains a stable orthogonal frame of format $(n,k)$
  with $n \geq 3$. It follows by Proposition \ref{p*ringembeddingrho}
  that $S$ is faithfully representable. 
\end{bew2}

\begin{prop}\lab{psimple*regringisrep} Every simple
  $*$-regular ring $R$ admits a faithful linear
  representation provided that $\LR$ does so.
\end{prop}

\begin{bew} Let $R$ be a simple $*$-regular ring $R$ without
  unit. Consider the set $P(R)$ of all projections in $R$. Since $R$
  is $*$-regular, $P(R)$ is a lattice. In particular, it is a directed 
set. By Lemma \ref{lprojx}, we have that $R$ is the directed union of
its subrings $R_e$, $e \in P(R)$. By Lemma
\ref{lultraprodofsubrings}, $R$ is a $*$-regular subring of an
ultraproduct of the $R_e$, $e \in P(R)$.
By Lemma \ref{leResimple}, for each projection $e$, the ring $R_e$ 
is a simple $*$-regular ring with unit $e$ and $\bar L(R_{eR_e})
 \cong [0,eR]$ is  representable.
 By Corollary
\ref{csimple*regringwithunitisrep}, each  $R_e$ is faithfully
representable. Hence, we can conclude with Lemma \ref{lsubringrep} and
Proposition \ref{pultraprodrep} that $R$ has a faithful
representation.  
\end{bew}

\subsubsection{Representations of Ideals}\lab{ssRepsofIdeals}

\begin{obs}Let $I$ be a two-sided ideal of the $*$-regular ring
  $R$. We can consider $I$ as a $*$-regular ring (without unit, if $I$
  is non-trivial) on its own; hence, we can consider representations of $I$.  
\end{obs}

\begin{prop}\lab{pidealrepringrep} Let $R$ be a $*$-regular ring and $I$ a
  two-sided ideal in 
  $R$ with a representation $ \varrho: I \ra End(V_D, \lan \cdot ,
  \cdot \ran)$. Denote the set of all projections in $I$ by $P(I)$, abbreviate 
$V_p:= \varrho(p)[V]$ for a projection $p \in P(I)$ and set
$$\rho(r):= \bigcup\limits_{p \in P(I)} \varrho(rp)_{|V_p}$$ 
Then $\rho$ is a representation of $R$ an appropriate
subspace $U$ of $V$, where the scalar product on $U_D$ is given by
restriction.
\end{prop}

\begin{bew} First, we have to show that the given definition of $\rho$
indeed defines a map $\rho: R \ra End(V_D)$. Recalling that the set of all
projections of a $*$-regular ring is directed, we consider two
projections $e,f \in P(I)$ with $e \leq f$, that is, with $e=fe$. We
have to show that the restrictions coincide on $V_e$, i.e., $\varrho(rf)_{|V_e}
= \varrho(re)_{|V_e}$. Since $e = fe$, we have 
$$\varrho(re)_{|V_e} = \varrho(rfe)_{|V_e} = (\varrho(rf) \circ
\varrho(e))_{|V_e} = \varrho(rf)_{|V_e},$$
as desired.

Second, we have to show that the map $\rho: R \ra End(V_D)$ is a
$*$-ring-homomorphism. For $0$ in $R$, we have
$$\rho(0) =  \bigcup\limits_{p \in P(I)} \varrho(0p)_{|V_p} = 0_V.$$ 
If $1 \in R$, we have 
$$\rho(1) =  \bigcup\limits_{p \in P(I)} \varrho(1p)_{|V_p} = 1_U$$ 
with $U:= \bigcup_{p \in P(I)} V_p$, that is, $\rho[R]$ acts on the subspace
$U$ of $V_D$.

For addition, let $r,s \in R$. We have
\begin{eqnarray*}
\rho(r+s) 
&=& \bigcup\limits_{p \in P(I)}\varrho((r+s)p)_{|V_p}= \bigcup\limits_{p
  \in P(I)} \varrho(rp + sp)_{|V_p}\\ 
&=& \bigcup\limits_{p \in P(I)} \varrho(rp)_{|V_p} + \varrho(sp)_{|V_p}
 = \bigcup\limits_{p \in P(I)} \varrho(rp)_{|V_p} + \bigcup\limits_{q
  \in P(I)} \varrho(sq)_{|V_q}\\ 
&=& \rho(r) + \rho(s).
\end{eqnarray*}
Therefore,  $\rho(r+s) = \rho(r) +
\rho(s)$ and $\rho(-r) = -\rho(r)$for all $r,s \in R$.

For multiplication, let $r,s$ be in $R$. We note that for each
$p \in P(I)$, there exists a $q_p \in P(I)$ such that $sp=q_psp$. We
claim that  
$$ \bigcup\limits_{p \in P(I)} (\varrho(rq_p) \circ \varrho(sp))_{|V_p}
 = \bigcup\limits_{q \in P(I)}
 \varrho(rq)_{|V_q} \circ  \bigcup\limits_{p \in P(I)}
 \varrho(sp)_{|V_p}.$$

\begin{bew2}Take $v \in U$. Then there exists $p_v \in P(I)$ with $v
 \in V_{p_v}$, so on the one hand
$$  \bigg(\bigcup\limits_{p \in P(I)} \big(\varrho(rq_p) \circ
\varrho(sp)\big)_{|V_p}\bigg) (v) = \big(\varrho(rq_{p_v})v \circ \varrho (s
p_v)\big) (v),$$ 
while on the other hand
\begin{eqnarray*}
&&\bigg(\bigcup\limits_{q \in P(I)} \varrho(rq)_{|V_q} \circ
 \bigcup\limits_{p \in P(I)} \varrho(sp)_{|V_p}\bigg)(v) =
 \bigcup\limits_{q \in P(I)} \varrho(rq)_{|V_q} \Big( \bigcup\limits_{p
 \in P(I)} \varrho(sp)_{|V_p}(v) \Big)\\ &=& 
\bigcup\limits_{q \in P(I)} \varrho(rq)_{|V_q}\Big(\varrho(sp_v)(v)\Big) =
 \varrho(rq_{p_v}) (\varrho (sp_v)(v)).
\end{eqnarray*}
\end{bew2}

Thus, we have 
\begin{eqnarray*}
\rho(r s) 
&=& \bigcup\limits_{p \in P(I)}\varrho((rs)p)_{|V_p}= \bigcup\limits_{p
  \in P(I)} \varrho(rq_psp)_{|V_p}\\ 
&=& \bigcup\limits_{p \in P(I)} (\varrho(rq_p) \circ \varrho(sp))_{|V_p}
 = \bigcup\limits_{q \in P(I)}
 \varrho(rq)_{|V_q} \circ  \bigcup\limits_{p \in P(I)} \varrho(sp)_{|V_p}\\ 
&=& \rho(r) \circ \rho(s).
\end{eqnarray*}\\

Now we examine the involution on $R$. For $r \in R$, consider $v,w \in
V$. Then take $e \in P(I)$ with $v 
\in V_e$. There exists $f_1 \in P(I)$ such that $f_1 re = re$ and $f_2
\in P(I)$ such that $w \in V_{f_2}$. Choosing $f: = f_1 \vee f_2$, we have 
\begin{eqnarray*}\lan \rho(r) v, w \ran  &=& 
\lan \bigcup\limits_{p \in P(I)} \varrho(rp)_{|V_p}v, w \ran = 
\lan \varrho(re) v, w \ran = \lan \varrho(fre) v, w \ran\\
&=&\lan v, \varrho(er^*f) w
\ran = \lan  v, \varrho(e) \varrho(r^*f) w \ran  = \lan
\varrho(e) v,  \varrho(r^*f) w \ran\\
&=& \lan v, \varrho(r^*f) w \ran   = \lan v, \rho(r^*)w \ran,
\end{eqnarray*}
where we have used that $\varrho: I \ra End(V_D)$ is a
$*$-ring-homomorphism and $v \in V_e, w \in V_f$. \end{bew}

\begin{lem}\lab{lidealrepringrep} Let $R$ be a $*$-regular ring and $I$ a
  two-sided ideal in $R$ with a representation $ \varrho: L \ra End(V_D, \lan
  .,. \ran)$. Denote the action of $R$ on the ideal $I$ given by left
  multiplication by $\lambda_I$, that is
$$\lambda_I: R \ra End(I_I) \qquad  \lambda_I(r) (x): = \widehat r (x)
= rx.$$
If the representation $\varrho: I \ra End(V_D, \lan \cdot , \cdot
\ran)$ is faithful and $\lambda_I: R \ra End(I_I)$ is injective, then the
representation $\rho: R  \ra End(V_D, \lan \cdot , \cdot \ran)$
defined in Proposition \ref{pidealrepringrep}  is faithful.
\end{lem}

\begin{bew2} Assume that $\rho(r) = 0$. This is equivalent to $\varrho(re) =
0$ for all $e \in P(I)$. As $\varrho: I \ra End(V_D, \lan \cdot ,
\cdot \ran)$ is faithful, this means that $re = 0$ for all $e \in
P(I)$. Since $I$ is a $*$-regular ring, for every element $x \in I$ there 
exists $e \in P(I)$ such that $ex = x$. Hence, we have that $rx = 0$
for all $x \in I$. Since we assumed the action of $R$ on $I$ given by left
multiplication to be injective, we have that $r=0$. This shows
that $\rho$ is injective.
\end{bew2}

\subsubsection{Subdirectly Irreducible $*$-Regular Rings}

In this section, we will show that each subdirectly irreducible
$*$-regular ring $R$ has a faithful representation
provided that $\LR$  does so.
\begin{obs} We may assume that $R$ is non-Artinian: Since every
  regular ring is semi-prime, a subdirectly irreducible $*$-regular
  ring which is Artinian is semi-simple, hence representable. 

Furthermore, the minimal two-sided ideal of $R$ is
non-Artinian, too (see \cite{ExistVarRRCML}, Proposition 2).

\end{obs}

\begin{prop}\lab{pactiononJ} Let $R$ be a subdirectly irreducible
  $*$-regular ring and let $J$ be the minimal two-sided ideal of
  $R$. Then the action $\lambda_J: R \ra End(J_J)$ of $R$ defined by  
$$\lambda_J(r)(x) = \hat{r}(x) = rx$$
is injective. 
\end{prop}

\begin{bew}Consider the left annihilator $A:=ann_R^l(J)$ of $J$ in $R$.
While a priori $A$ is only a left ideal, it can be shown that $A$ is
indeed closed under left and right multiplication by
elements of $R$. Since $A$ is closed under addition, $A$ is a
two-sided ideal in $R$. Since $J$ does not annihilate itself,  we can conclude
that $A$ is trivial. Therefore, the action of $R$ on $J$ defined by
left multiplication is injective.
\end{bew}

\begin{lem}\lab{lJsimple} The minimal ideal $J$ of a subdirectly irreducible
  $*$-regular ring $R$ is a simple $*$-regular ring.
\end{lem}

\begin{bew2} For a non-vanishing ideal $A$ in $J$, consider the ideal
  generated by $A$ in $R$. 
\end{bew2}

\begin{rem}Note that one does not need that $R$ contains a unit.
\end{rem}

\begin{lem}\lab{leIResimple} Let $R$ be a $*$-regular ring and $I$ a
  minimal two-sided ideal in 
  $R$, in particular simple as a ring. Let $e$ be a projection in
  $I$. 

Then the ring $R_e = eRe$ is simple. 
\end{lem}

\begin{bew2} For a non-vanishing ideal $A$ in $R_e$, consider the
  ideal generated by $A$ in $I$. 
\end{bew2}

\begin{prop}\lab{pJrep} Considered as simple $*$-regular ring, the
  minimal ideal $J$ of a
 subdirectly irreducible $*$-regular ring $R$
  has a faithful representation provided that $\bar L(I_R)$ 
does so.  
\end{prop}

\begin{bew2} Proposition \ref{psimple*regringisrep}.
\end{bew2}

\begin{thm}\lab{csi*regringrep}Every  subdirectly irreducible
  $*$-regular ring $R$ has a faithful representation
provided that $\LR$ does so.
\end{thm}
\begin{bew2}Lemma \ref{lidealrepringrep}.
\end{bew2}


Authors address: Fachbereich Mathematik der Technischen Universit\"at Darmstadt
\end{document}